\documentclass[11pt]{article}
\usepackage {pstcol}
\usepackage {pstricks,pst-node}    

\usepackage{mathtools}
\usepackage{amsmath,amssymb,amsbsy}
\usepackage{graphicx}
\usepackage{multirow}

\usepackage{amsfonts}
\usepackage{amsthm}
\usepackage{amssymb}
\usepackage{graphicx, enumerate}
\usepackage{amsmath}
\usepackage{latexsym}
\usepackage{longtable}
\usepackage{tabularx}
\usepackage{amsmath}
\usepackage{amsfonts}
\usepackage{amsthm}
\usepackage{setspace}
\usepackage{graphicx}
\usepackage{float}
\usepackage{rotating}
\usepackage{tikz}
\usepackage{verbatim}
\usepackage{soul}

\def\zet{\mathbb{Z}}

\def\fk2{\lfloor\frac{k}{2}\rfloor}
\def\ck2{\lceil\frac{k}{2}\rceil}
\def\k1{\lfloor\frac{k+1}{2}\rfloor}

\def\deg {\mathrm{deg}}

\newtheorem{theorem}{Theorem}[section]
\newtheorem{lemma}[theorem]{Lemma}
\newtheorem{conjecture}[theorem]{Conjecture}
\newtheorem{proposition}[theorem]{Proposition}

\newcommand{\trou}{\vspace{5 mm} \noindent}

\newcommand{\cbdo}{\hfill \rule{.1in}{.1in}}    
\newcommand{\prf}{\noindent{\bf Proof.}\ }

\title{On rainbow caterpillars in elementary $p$-groups}
\author{Sylwia Cichacz, Barbara Krupińska, Mariusz Woźniak\\\\
\normalsize AGH University, \vspace{2mm} Cracow, Poland}

\begin{document}

\maketitle
\begin{abstract} Given a finite Abelian group $(A,+)$, consider a tree $T$ with $|A|$ vertices. The labeling $f \colon V (T) \rightarrow A$ of the vertices of some graph $G$ induces an edge labeling in $G$, thus the edge $uv$ receives the label $f (u) + f (v)$. 
 The tree $T$ is $A$-rainbow colored if $f$ is a bijection and edges have different colors. In this paper, we give necessary and sufficient conditions for a caterpillar with three spine vertices to be $A$-rainbow, when $A$ is an elementary $p$-group.
\end{abstract}

\section{Introduction}

Given two integers $a$, $b$, we denote by $[a,b]$  the set consisting of the integers $a,a+1,\ldots,b$. Recall that the exponent of a group $A$ ($\exp(A)$) is defined as the least common multiple of the orders of all elements of the group.  An elementary Abelian group is an Abelian group in which all elements other than the identity have the same order. This common order must be a prime number. Thus, we call $A$
 an \textit{elementary $p$-group} if and only if $\exp(A)=p$, i.e. $G\cong (\zet_p)^k$ for some prime number $p$. By $\langle a_1,\ldots,a_s \rangle$ we denote a subgroup generated by $ a_1,\ldots,a_s$ in the group $A$. Let $\langle S\rangle=\langle a_1,\ldots,a_s \rangle$ for $S=\{a_1,\ldots,a_s\}$.

 The labeling $f \colon V (T) \rightarrow A$ of the vertices of some graph $G$ induces an edge labeling in $G$, thus the edge $uv$ receives the label $f (u) + f (v)$. 
 The tree $T$ in $A$-rainbow colored if $f$ is a bijection and edges have different colors.

In \cite{JamKin} Jamison and Kinnersley considered an $A$-rainbow labeling of trees. They proved the following:
\begin{theorem}[\cite{JamKin}]
Let $T$ be a tree of order $n$ and $A$ be an Abelian group such that $|A|=n$ and $\exp(A)=m>2$.
If $T$ has adjacent vertices $u$ and $v$ such that
\begin{itemize}
    \item  $\deg(u) \equiv \deg(v) \equiv 0 \pmod m$,
\item  $\deg(x) \equiv 1 \pmod m$ for all $x \in V(T ) \setminus \{u, v\}$, and
\item $uv \in E(T)$,
\end{itemize}
then $T$ is not $A$-rainbow.
\end{theorem}

Jamison and Kinnersley also studied caterpillars. A \textit{caterpillar} is either $K_2$ or a tree on at least 3 vertices such that deleting its leaves we obtain a path of order at least $1$.
Let $C(h_1,\ldots, h_s)$ denote the caterpillar
with $s$ spine vertices and with $h_i$ hairs on the $i$th spine vertex.

\begin{theorem}[\cite{JamKin}]\label{stonogi1}
Let $T\cong C(h_1,h_2,\ldots,h_s)$ be a caterpillar of order $n$ and $A$ be an Abelian group such that $|A|=n$ and $\exp(A)=m>2$. Let $r_1,\ldots , r_s$ be the residues of $h_1, \ldots , h_s$ modulo $m$. If $r_1+\ldots + r_s+s=m$, then $T$ has an $A$-rainbow labeling.

\end{theorem}
\begin{theorem}[\cite{JamKin}]\label{stonogi2}
Let $T$ be a caterpillar of order $n$ and $A$ be an Abelian group such that $|A|=n$ and $\exp(A)=m>2$.
\begin{itemize}
    \item $T\cong C(h_1, h_2)$ has an $A$-rainbow labeling if and only if $h_1 \not\equiv -1 \pmod m.$
    \item $T\cong C(h_1,0, h_3)$ has an $A$-rainbow labeling if and only if $h_1 \not\equiv -1 \pmod m$ and $h_3 \not\equiv -1 \pmod m.$
\end{itemize}
\end{theorem}

Many related types of labelings have been studied in the literature. The most connected to the $A$-rainbow is $A$-cordial labeling, which was introduced by Hovey as a generalization of cordial and harmonious labeling \cite{Hovey}. If $A$ is an Abelian group, then the labeling $f \colon V (G) \rightarrow A$ of the vertices of some graph $G$ induces an edge labeling in $G$, thus the edge $uv$ receives the label $f (u) + f (v)$. A graph $G$ is $A$-cordial if there is a vertex-labeling such that (1) the vertex label classes differ in size by at most one and (2) the induced edge label classes differ in size by at most one.  Note that $A$-rainbow labeling is $A$-cordial for a tree with $|A|$-vertices.

In the literature, cordial labeling in cyclic groups is the main focus of studies. There is a famous (still open) conjecture that states that all trees are $\zet_k$-cordial
for all $k$ \cite{Hovey}.  The conjecture is confirmed for caterpillars:
\begin{theorem}[\cite{Hovey}]\label{kcordial}
  Caterpillars are $\mathbb{Z}_k$-cordial for all $k > 0$.
\end{theorem}

 The situation changes a lot if $A$ is not cyclic. It was proven that all trees, except of $P_4$ and $P_5$, are $\mathbb{Z}_2^2$-cordial \cite{Klein,Pechenik2}.

A family of $\mathbb{G}$ of graphs is called weakly \textit{$A$-cordial} if all but a finite number of elements of  $\mathbb{G}$ are $A$-cordial.  Erickson et al. posted the conjecture \cite{Klein}:
 \begin{conjecture}[\cite{Klein}]\label{Klein}For any Abelian group $A$, the set of trees is weakly $A$-cordial.
 \end{conjecture}
Observe that Conjecture \ref{Klein} is not true for general trees, even if we consider $A$-rainbow coloring instead of $A$-cordial labeling (i.e., $A$-cordial labeling in which $|A|=|V(G)|$). Namely, let $m>2$ be an integer, $A\cong (\zet_{m^2})^k \times(\zet_m)^l$, and $T \cong C(a,b)$ for $a \equiv-1 \pmod{m^2}$  or $T \cong C(a,0,b)$, for $a,b \equiv-1 \pmod{m^2}$ then based on Theorem~\ref{stonogi2}  we can clearly see that this set of trees is not weakly $A$-cordial.
 
 In this paper, we give necessary and sufficient conditions for a caterpillar with three spine vertices to be $A$-rainbow, when $A$ is an elementary $p$-group. As a result, we present a new family of trees that serve as counterexamples to Conjecture \ref{Klein}.


\section{Extension of a labeling}


\subsection{Notation and formulation of the main result}
 Suppose that we have labeled the spine of  a caterpillar $C( |X_1|, \ldots, |X_s|)$ by  $a_1,\ldots, a_s\in A$. Let now $C[a_1,\ldots, a_s; X_1, \ldots, X_s]$ denote such caterpillar with $s$ spine vertices where $a_i\in A$ is a label given to $i$-th vertex of the spine and $X_i$ is a set of hairs connected to this vertex.
 We will extend the labeling given to the spine (i.e. we will label vertices from $X_1\cup\ldots\cup X_s$) to an $A$-rainbow labeling of the caterpillar $C( |X_1|, \ldots, |X_s|)$. To do that we need the definition of Cayley digraphs, which  are defined with a group $A$  and a subset $S$ 
 of $A$, the vertices of the Cayley digraph $Cay(A,S)$
 are the elements of the group, and its arcs are all the couples $(a,as)$
 with $a\in A$  and $s\in S$ \cite{Del}. Observe that $Cay(A,S)$ is weakly connected if and only if
$\langle S \rangle=A$. For $\langle S \rangle\neq A$ the weakly connected components correspond to the cosets of $\langle S \rangle$.

The purpose of this section is to provide a complete characterization of those caterpillars with a 3-vertex spine that have a rainbow coloring for the group  $A = {\zet}_p^k$, where $p$ is prime.

So, this is a tree  $T\cong C[a_1,a_2,a_3;X,Y,Z]$. 
Let $\alpha=|X| \pmod p$, $\beta=|Y| \pmod p$ and $\gamma=|Z| \pmod p$, where $\alpha,\beta,\gamma\in\{0,1,\ldots,p-1\}$. 
Let us note that $\alpha + \beta + \gamma = (p-3)\pmod p$.

The existence of a rainbow coloring depends on the triple $(\alpha, \beta, \gamma)$ and on the labels $[a_1,a_2,a_3]$ of the spine vertices. 

In this section, we will use a somewhat simplified nomenclature.
Namely, if there exists an $A$-rainbow coloring of the tree $T\cong C[a_1,a_2,a_3;X,Y,Z]$, then we will say that the triple $(\alpha, \beta, \gamma)$ is realizable in the model $[a_1,a_2,a_3]$.

\begin{figure}[ht]
\begin{center}
\includegraphics[width=16cm]{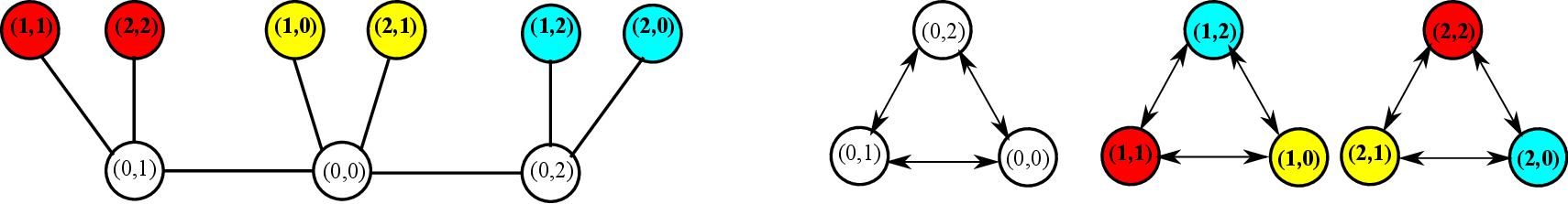}
\caption{A $\mathbb{Z}_3^2$-rainbow labeling of  $C[(0,1),(0,0),(0,2); X, Y, Z]$ and a corresponding partition of $Cay(\mathbb{Z}_3^2,\{(0,1),(0,0),(0,2)\})$}\label{c12}
\end{center}
\end{figure}

Note that in order to define a rainbow coloring, it is enough to indicate for any element of the group $A$, i.e. for any vertex in the graph $\mbox{Cay}(A,S)$ (except for the set $\{a_1,a_2,a_3\}$), to which of the sets $X, Y, Z$ it should belong (see Figure~\ref{c12}). { Therefore, assigning a label $x$ ($y$, $z$, respectively) to  $u\in V(\mbox{Cay}(A,S))$ means that in the partition of $\mbox{Cay}(A,S)$ the vertex $u$ is the image of some vertex of the tree $T$, which is a leaf from the set $X$ ($Y$, $Z$, respectively).}

The component of the graph $\mbox{Cay}(A,S)$,
where $S=\{a_1,a_2,a_3\}$, containing the vertices $\{a_1,a_2,a_3\}$ will be called the \emph{spine} component, and the remaining components will be called the\emph{ regular} components. 

\begin{lemma}\label{sym}
If the triple $(\alpha, \beta, \gamma)$ is realizable, then 
the triple $(\gamma, \beta, \alpha)$ is also realizable.
\end{lemma}

\prf Just change the model $[a_1,a_2,a_3]$ to $[a_3,a_2,a_1]$. \cbdo

\begin{lemma}\label{zero}
If the triple $(\alpha, \beta, \gamma)$ is realizable in the model $[a_1,a_2,a_3]$, then it is also realizable in the model $[a,0,b]$, for some $a$ and $b$.
\end{lemma}

\prf Note that adding a fixed element $A$ to all vertices of the graph $\mbox{Cay}(A,S)$ changes the value of the sums on the edges, but the different sums remain different. Adding an element $-a_2$ to all vertices gives us the model $[a_1-a_2,0,a_3-a_2]$.  \cbdo

\trou

We will continue to consider the model $[a,0,b]$. 
Note that the vertex $u=b-a$ cannot be labeled by $x$, since then on edge $ua$ we would have the same sum as on edge $0b$. We express this fact by saying that the assignment $x=b-a$ is forbidden. In turn, the symbol $x \stackrel{a}{\longrightarrow} y$ means that if a vertex $u$ gets label $x$, then the vertex $u+a$ cannot get the label $y$. 
The following lemma presents the remaining forbidden assignments, the easy proof of which we leave to the reader.

\begin{lemma}\label{FA}
In the model $[a,0,b]$, we have the following forbidden assignments (FA).

FA: $x=b-a$, $z=a-b$; $x \stackrel{a}{\longrightarrow} y$, $z \stackrel{b}{\longrightarrow} y$, 
$z \stackrel{b-a}{\longrightarrow} x$.

\cbdo
\end{lemma}

\begin{lemma}\label{xyz}
In a regular component, in the model $[a,0,b]$ where $b\in  \langle a\rangle$, either all three elements $x$, $y$ {and $z$} are present, or only one of them.
\end{lemma}

\prf The proof follows from the fact that $p$ is a prime number.  \cbdo

\trou
Most often, we will deal with the model $[a,0,b]$ where $b\in \langle a\rangle$. However, we will need the result below, which refers to the situation when $b\notin  \langle a\rangle$. In such a model, each component has $p^2$ vertices. The spine component contains $0$ and is generated by elements $a$ and $b$. The vertices are therefore of the form $sa + tb$ {where $s,t \in [0,p-1]$}

\begin{lemma}\label{krata}
Let $\alpha + \beta + \gamma \equiv (p-3)\pmod p$. If the triple $(\alpha, \beta, \gamma)$ is realizable  in the model $[a,0,b]$ where $b\not\in \langle a\rangle$, then $|Y| \ge p-1$.

\end{lemma}

\prf  
Consider an oriented path of length $p-2$ starting at the vertex $a-b$ generated by $a-b$ of the form
$(a-b)+s(a-b)$, $0\le s \le p-2$. Note that the last vertex on this path is
$(a-b)+(p-2)(a-b)=b-a$.

We will show that at least one vertex on this path has the label $y$. Suppose all of them have labels $x$ or $z$.
Note that the vertex $a-b$ must have label $x$ because $z$ is forbidden for this vertex. 
Similarly, vertex $b-a$ must have label $z$. But then there must be a vertex with the label $x$ such that the next vertex on the path has the label $z$.
But then we have the situation $x \stackrel{a-b}{\longrightarrow} z$, which is forbidden.

Therefore, there exists such $s$ that a vertex $w=(a-b)+s(a-b)$ has the label $y$.

Due to forbidden assignments, the vertex $w-a$ cannot have the label $x$, and the vertex $w-b$ cannot have the label $y$. Note that both vertices belong to the directed cycle generated by $a-b$. This cycle, therefore, contains vertices with both labels $x$ and $z$. If this cycle contained only vertices with such labels, we would have a forbidden situation.

It must therefore contain at least one vertex with the label $y$.

Continuing this reasoning, we obtain that there are at least $p-1$ vertices with label $y$.\cbdo

\begin{figure}
\psset{unit=1cm}
\psset{radius=0.2}

\begin{pspicture}(13.5,6)

\put(3.5,0.5)
{\begin{pspicture}(0,0)

\rput(0,6){\dianode[linestyle=none,fillstyle=solid,fillcolor=yellow]{C11}{\mbox{0}}}
\rput(1.5,6){\dianode[linestyle=none,fillstyle=solid,fillcolor=yellow]{C12}{\mbox{a}}}
\rput(0,4.5){\dianode[linestyle=none,fillstyle=solid,fillcolor=yellow]{C21}{\mbox{b}}}

\rput(3,6){\circlenode[linestyle=none,fillstyle=solid,fillcolor=green]{C13}{\hspace{2mm}}}
\rput(4.5,6){\circlenode[linestyle=none,fillstyle=solid,fillcolor=green]{C14}{\hspace{2mm}}}
\rput(6,6){\circlenode[linestyle=none,fillstyle=solid,fillcolor=green]{C15}{\hspace{2mm}}}

\rput(1.5,4.5){\circlenode[linestyle=none,fillstyle=solid,fillcolor=green]{C22}{\hspace{2mm}}}
\rput(3,4.5){\circlenode[linestyle=none,fillstyle=solid,fillcolor=green]{C23}{\hspace{2mm}}}
\rput(4.5,4.5){\circlenode[linestyle=none,fillstyle=solid,fillcolor=green]{C24}{\hspace{2mm}}}
\rput(6,4.5){\circlenode[linestyle=none,fillstyle=solid,fillcolor=green]{C25}{\mbox{z}}}

\rput(0,3){\circlenode[linestyle=none,fillstyle=solid,fillcolor=green]{C31}{\hspace{2mm}}}
\rput(1.5,3){\circlenode[linestyle=none,fillstyle=solid,fillcolor=green]{C32}{\hspace{2mm}}}
\rput(3,3){\circlenode[linestyle=none,fillstyle=solid,fillcolor=green]{C33}{\mbox{x}}}
\rput(4.5,3){\circlenode[linestyle=none,fillstyle=solid,fillcolor=green]{C34}{\hspace{2mm}}}
\rput(6,3){\circlenode[linestyle=none,fillstyle=solid,fillcolor=green]{C35}{\hspace{2mm}}}

\rput(0,1.5){\circlenode[linestyle=none,fillstyle=solid,fillcolor=green]{C41}{\hspace{2mm}}}
\rput(1.5,1.5){\circlenode[linestyle=none,fillstyle=solid,fillcolor=green]{C42}{\mbox{z}}}
\rput(3,1.5){\circlenode[linestyle=none,fillstyle=solid,fillcolor=red]{C43}{\mbox{y}}}
\rput(4.5,1.5){\circlenode[linestyle=none,fillstyle=solid,fillcolor=green]{C44}{\hspace{2mm}}}
\rput(6,1.5){\circlenode[linestyle=none,fillstyle=solid,fillcolor=green]{C45}{\hspace{2mm}}}

\rput(0,0){\circlenode[linestyle=none,fillstyle=solid,fillcolor=green]{C51}{\hspace{2mm}}}
\rput(1.5,0){\circlenode[linestyle=none,fillstyle=solid,fillcolor=green]{C52}{\mbox{x}}}
\rput(3,0){\circlenode[linestyle=none,fillstyle=solid,fillcolor=green]{C53}{\hspace{2mm}}}
\rput(4.5,0){\circlenode[linestyle=none,fillstyle=solid,fillcolor=green]{C54}{\hspace{2mm}}}
\rput(6,0){\circlenode[linestyle=none,fillstyle=solid,fillcolor=green]{C55}{\hspace{2mm}}}

\ncarc{->}{C11}{C12}
\ncarc{->}{C12}{C13}
\ncarc{->}{C13}{C14}
\ncarc{->}{C14}{C15}
\ncarc{->}{C15}{C11}

\ncarc{->}{C21}{C22}
\ncarc{->}{C22}{C23}
\ncarc{->}{C23}{C24}
\ncarc{->}{C24}{C25}
\ncarc{->}{C25}{C21}

\ncarc{->}{C31}{C32}
\ncarc{->}{C32}{C33}
\ncarc{->}{C33}{C34}
\ncarc{->}{C34}{C35}
\ncarc{->}{C35}{C31}

\ncarc{->}{C41}{C42}
\ncarc{->}{C42}{C43}
\ncarc{->}{C43}{C44}
\ncarc{->}{C44}{C45}
\ncarc{->}{C45}{C41}

\ncarc{->}{C51}{C52}
\ncarc{->}{C52}{C53}
\ncarc{->}{C53}{C54}
\ncarc{->}{C54}{C55}
\ncarc{->}{C55}{C51}

\ncarc{->}{C11}{C21}
\ncarc{->}{C21}{C31}
\ncarc{->}{C31}{C41}
\ncarc{->}{C41}{C51}
\ncarc{->}{C51}{C11}

\ncarc{->}{C12}{C22}
\ncarc{->}{C22}{C32}
\ncarc{->}{C32}{C42}
\ncarc{->}{C42}{C52}
\ncarc{->}{C52}{C12}

\ncarc{->}{C13}{C23}
\ncarc{->}{C23}{C33}
\ncarc{->}{C33}{C43}
\ncarc{->}{C43}{C53}
\ncarc{->}{C53}{C13}

\ncarc{->}{C14}{C24}
\ncarc{->}{C24}{C34}
\ncarc{->}{C34}{C44}
\ncarc{->}{C44}{C54}
\ncarc{->}{C54}{C14}

\ncarc{->}{C15}{C25}
\ncarc{->}{C25}{C35}
\ncarc{->}{C35}{C45}
\ncarc{->}{C45}{C55}
\ncarc{->}{C55}{C15}

\ncarc[linestyle=dashed]{->}{C52}{C43}
\ncarc[linestyle=dashed]{->}{C43}{C34}
\ncarc[linestyle=dashed]{->}{C34}{C25}

\ncarc[linestyle=dotted]{->}{C51}{C42}
\ncarc[linestyle=dotted]{->}{C42}{C33}
\ncarc[linestyle=dotted]{->}{C33}{C24}
\ncarc[linestyle=dotted]{->}{C24}{C15}
\ncarc[linestyle=dotted]{->}{C15}{C51}

\put(1.1,-0.5){$a-b$}
\put(6.3,4.4){$b-a$}

\end{pspicture}}
\end{pspicture}
\caption {Spine component in the model $[a,0,b]$ where $b \notin \langle a \rangle$} \label{5x5}
\end{figure}

\begin{theorem}\label{p5}
If $T$ is a caterpillar with a spine of length two and order $|V(T)|=p^k$, $k>1$ and $A = {\zet}_p^k$, where $p$ is prime greater than or equal to 5. Ta hen $T$ is $A$-rainbow except for the following cases:

\begin{enumerate}
    \item  $\beta=p-2$, and $\alpha=0$, $\gamma=p-1$, or, symmetrically,
$\alpha=p-1$, $\gamma=0$, 

\item $|Y|=0$, and $\alpha=p-1$, $\gamma=p-2$, or, symmetrically,
$\alpha=p-2$, $\gamma=p-1$,

\item $|Y|=1$, and $\alpha=p-1$, $\gamma=p-3$, or, symmetrically,
$\alpha=p-3$, $\gamma=p-1$.\end{enumerate}

\end{theorem}

The proof of the theorem is divided into several parts and will be presented in the next 
subsections.

\subsection{Case $\alpha=\gamma$}\label{xeqy}
\begin{figure}
\psset{unit=0.96cm}
\psset{radius=0.2}

\begin{pspicture}(12,2)

\put(0.5,0.5)
{\begin{pspicture}(0,0)

\put(3,0.7){$\overbrace{\hspace{3cm}}$}
\put(4,1){$\beta~ \mbox{times}$}
\put(4.3,0){$\cdots$}

\put(7.5,0.7){$\overbrace{\hspace{4.5cm}}$}
\put(7.5,0.3){$\overbrace{\hspace{1.5cm}}$}
\put(10.5,0.3){$\overbrace{\hspace{1.5cm}}$}

\put(9.3,1){$\alpha~ \mbox{times}$}
\put(9.5,0){$\cdots$}

\ncarc{->}{X5}{X6}

\rput(0,0){\dianode[linestyle=none,fillstyle=solid,fillcolor=yellow]{C1}{\mbox{0}}}
\rput(1.5,0){\dianode[linestyle=none,fillstyle=solid,fillcolor=yellow]{C2}{\mbox{a}}}
\rput(13.5,0){\dianode[linestyle=none,fillstyle=solid,fillcolor=yellow]{C10}{\mbox{-a}}}

\rput(3,0){\circlenode[linestyle=none,fillstyle=solid,fillcolor=green]{C3}{\mbox{y}}}
\rput(6,0){\circlenode[linestyle=none,fillstyle=solid,fillcolor=green]{C5}{\mbox{y}}}

\rput(7.5,0){\circlenode[linestyle=none,fillstyle=solid,fillcolor=green]{C6}{\mbox{x}}}
\rput(9,0){\circlenode[linestyle=none,fillstyle=solid,fillcolor=green]{C7}{\mbox{z}}}
\rput(10.5,0){\circlenode[linestyle=none,fillstyle=solid,fillcolor=green]{C8}{\mbox{x}}}
\rput(12,0){\circlenode[linestyle=none,fillstyle=solid,fillcolor=green]{C9}{\mbox{z}}}

\ncarc{->}{C1}{C2}
\ncarc{->}{C2}{C3}
\ncarc{->}{C5}{C6}
\ncarc{->}{C6}{C7}
\ncarc{->}{C8}{C9}
\ncarc{->}{C9}{C10}
\ncarc{->}{C10}{C1}

\end{pspicture}}
\end{pspicture}
\caption {One of the possible assignments for the spine component in the model $[a,0,-a]$} \label{a0-a}
\end{figure}

In this case we use the $[a,0,-a]$ model. Then the forbidden assignments  are of the form:

\noindent FA: $x=-2a$, $z=2a$; $x \stackrel{a}{\longrightarrow} y$, $z \stackrel{-a}{\longrightarrow} y$, 
$z \stackrel{-2a}{\longrightarrow} x$.

We can also write the last two conditions in the form $y \stackrel{a}{\longrightarrow} z$,
$x \stackrel{2a}{\longrightarrow} z$. 

If for the triple $(\alpha,\beta,\alpha)$ there is $2\alpha + \beta = p-3$, 
then it can be realized on the spine component. One possible realization is shown in 
Fig.~\ref{a0-a}.

Now let's consider the case where $2\alpha + \beta = 2p-3$.  We are therefore dealing with a triple that can be written in the following form
$\left(\frac{2p-3-\beta}{2}, \beta, \frac{2p-3-\beta}{2}\right)$ where $\beta$ is an odd number less than or 
equal to $p-2$. 

Then, for the implementation, we need to use both the spine component and one regular component. 
The key to what follows is the observation that the triple $\left(\frac{p-1}{2}, 1, \frac{p-1}{2}\right)$ can be implemented on a regular component. Indeed, it is enough to place one $y$ anywhere on the cycle generated by element $a$, \emph{i.e.} on the cycle $u, u+a, u+ 2a, \ldots, u +(p-1)a$ and complete the cycle with $\frac{p-1}{2}$ pairs $xz$ in accordance with the cycle orientation.

Since the triple of the form
$\left(\frac{p-2-\beta}{2}, \beta{-1}, \frac{p-2-\beta}{2}\right)$ can be realized on the spine component, the possibility of realizing $\left(\frac{2p-3-\beta}{2}, \beta, \frac{2p-3-\beta}{2}\right)$ results from the following equality

$$\left(\frac{2p-3-\beta}{2}, \beta, \frac{2p-3-\beta}{2}\right)=\left(\frac{p-2-\beta}{2}, \beta{-1 }, \frac{p-2-\beta}{2}\right)+ \left(\frac{p-1}{2}, 1, \frac{p-1}{2}\right).$$

Finally, to realize the triple $(p-1,p-1,p-1)$ we need a spine component and two regular components. The possibility of realizing this triple follows from the equality

$$(p-1,p-1,p-1)=(0,p-3,0)+ 2\left(\frac{p-1}{2}, 1, \frac{p-1}{2}\right).$$

\subsection{Case $\beta \ge \gamma$ or $\beta \ge \alpha$}

In this case, we use the $[a,0,2a]$ model. Then the forbidden assignments  are of the form:

\noindent FA:  $z=2a$; $x \stackrel{a}{\longrightarrow} y$, $z \stackrel{2a}{\longrightarrow} y$, 
$z \stackrel{a}{\longrightarrow} x$.

\noindent{\bf Remark.} Note that the forbidden assignment for vertex i.e., $x=a$ 
is satisfied due to the choice of labels for the spine vertices.

This time let's start with a regular component in case there are 3 labels { $x$, $y$, and $z$.} 
We will look at this component as an oriented cycle generated by the element $a$. If we put $z$ somewhere, then the next element cannot be $x$ (since $z \stackrel{a}{\longrightarrow} x$ is forbidden). If the next element is $z$, then labeling the next vertex will not be possible, because we would have $z \stackrel{2a}{\longrightarrow} y$. So the next vertex must be labeled $y$. The next vertex cannot be $y$ for the same reasons as before. To sum up, if some vertex has label $z$, 
then the next one must have label $y$, and after a pair $zy$ comes another pair $zy$ or $x$.
In particular, it follows that on such a component it is possible to realize only such triples 
$(\alpha, \beta, \gamma)$ where $\beta=\gamma$. It also follows from the above that every such triple $(p- 2j, j,j)$ for $j\in\left\{0,1,\ldots, \frac{p-1}{2}\right\}$  can be realized on a regular component.

Let's put $r = \beta - \gamma$.
If the triple we have to realize is of the form $(\alpha,\gamma + r,\gamma)$ where $\alpha + 2\gamma + r = p-3$, then it can be realized on the spine component. One possible realization is shown in 
Fig.~\ref{a02a}.

Now let's consider the case where $\alpha + 2\gamma + r  = 2p-3$.  We are therefore dealing with a triple that can be written in the following form
$(2p-3-r - 2\gamma, \gamma + r, \gamma)$. In this case, we need to use both the spine component and one regular component.

 Our further action depends on $\gamma$.  Assume first $0 \le \gamma \le \frac{p-1}{2}$, then the  triple $(p- 2\gamma, \gamma, \gamma)$ can be realized on a regular component. Then, it is enough to realize 
$(p-3-r, r,0)$ on the spine component, and we are done. This is possible if $r\le p-3$. Note that if $r\geq p-2$, then $r\in\{p-2,p-1\}$, which implies $(\alpha,\beta,\gamma)\in\{(p-3,p-1,1),(p-1,p-2,0),(p-2,p-1,0)\}$. These cases we will consider later.

Assume now $\gamma \ge \frac{p+1}{2}$. {Note that $r=\beta-\gamma\leq p-3$.} The triple $(1,\frac{p-1}{2}, \frac{p-1}{2} )$ we realize on a regular component. If $\alpha>0$, then the triple $\left(p-4-\beta-\gamma,\beta-\frac{p-1}{2},\gamma- \frac{p-1}{2} \right)$ we realize on the spine component, since $\beta+\gamma\leq p-4$ in this case.  For $\alpha=0$, there is $\beta=p-1$ and $\gamma=p-2$. Recall that $|X|>0$.  Therefore, we could write this triple in the form $(p, p-1, p-2)$. In particular, it is clear that its realization requires three components, and the possibility of realization follows from the equality below.
 
$$(p, p-1, p-2) = (p-4, 1, 0) + \left(1, \frac{p-1}{2},  \frac{p-1}{2}\right) +\left(3, \frac{p-3}{2},  \frac{p-3}{2}\right).$$


If $r = p-2$ and $\gamma =1$  we are dealing with the triple $(p-3, p-1, 1)$. 
Then, instead of this triple, we consider the triple $(1, p-1, p-3)$ and we proceed as above in the case where $\gamma \ge \frac{p+1}{2}$.

Finally, let's consider the case when $\gamma =0$. There are two cases where the methods described above cannot be applied directly. Namely, it is about the triples 
$(p-1, p-2, 0)$ and $(p-2, p-1, 0)$. In the first case, rainbow labeling does not exist, 
as we will show later (cf. Lemma~\ref{x=0}). In the last case, again using Lemma~\ref{sym}, we consider the triple $(0, p-1, p-2)$ and we are done as above.

The last triple to consider is $(p-1,p-1,p-1)$, which was already considered previously. So we get another possibility of realization. To realize this triple we need a spine component and two regular components. The possibility of realizing this triple in the model $[a,0,2a]$ follows from the equality

$$(p-1,p-1,p-1)=(p-3,0,0)+ 2( 1, \frac{p-1}{2}, \frac{p-1}{2}).$$

\begin{figure}
\psset{unit=1.08cm}
\psset{radius=0.2}

\begin{pspicture}(13.5,1.5)

\put(0.5,0.5)
{\begin{pspicture}(0,0)

\put(3,0.7){$\overbrace{\hspace{1.5cm}}$}
\put(3,1){$r ~ \mbox{times}$}
\put(3.5,0){$\ldots$}

\put(5.5,0.7){$\overbrace{\hspace{1.5cm}}$}
\put(5.5,1){$\alpha ~ \mbox{times}$}
\put(6,0){$\ldots$}

\put(8,0.7){$\overbrace{\hspace{4cm}}$}
\put(8,0.3){$\overbrace{\hspace{1cm}}$}
\put(11,0.3){$\overbrace{\hspace{1cm}}$}

\put(9.3,1){$\gamma~ \mbox{times}$}
\put(9.7,0){$\cdots$}

\rput(0,0){\dianode[linestyle=none,fillstyle=solid,fillcolor=yellow]{C1}{\mbox{0}}}
\rput(1,0){\dianode[linestyle=none,fillstyle=solid,fillcolor=yellow]{C2}{\mbox{a}}}
\rput(2,0){\dianode[linestyle=none,fillstyle=solid,fillcolor=yellow]{C3}{\mbox{2a}}}

\rput(3,0){\circlenode[linestyle=none,fillstyle=solid,fillcolor=green]{C4}{\mbox{y}}}
\rput(4.5,0){\circlenode[linestyle=none,fillstyle=solid,fillcolor=green]{C5}{\mbox{y}}}

\rput(5.5,0){\circlenode[linestyle=none,fillstyle=solid,fillcolor=green]{C6}{\mbox{x}}}
\rput(7,0){\circlenode[linestyle=none,fillstyle=solid,fillcolor=green]{C7}{\mbox{x}}}

\rput(8,0){\circlenode[linestyle=none,fillstyle=solid,fillcolor=green]{C8}{\mbox{z}}}
\rput(9,0){\circlenode[linestyle=none,fillstyle=solid,fillcolor=green]{C9}{\mbox{y}}}

\rput(11,0){\circlenode[linestyle=none,fillstyle=solid,fillcolor=green]{C10}{\mbox{z}}}
\rput(12,0){\circlenode[linestyle=none,fillstyle=solid,fillcolor=green]{C11}{\mbox{y}}}

\ncarc{->}{C1}{C2}
\ncarc{->}{C2}{C3}
\ncarc{->}{C3}{C4}
\ncarc{->}{C5}{C6}
\ncarc{->}{C7}{C8}
\ncarc{->}{C8}{C9}
\ncarc{->}{C10}{C11}
\ncarc{->}{C11}{C1}

\end{pspicture}}
\end{pspicture}
\caption {One of the possible assignments for the spine component in the model $[a,0,2a]$} 
\label{a02a}
\end{figure}

\subsection{Case $\beta < \gamma$ and $\beta < \alpha$}

Unlike previous cases when one model was used to implement many triples, this time we will select an appropriate model for a given triple. 

Let $i\in A, i\neq 0$. The components of the graph generated by $i$ will be considered as oriented cycles of length $p$ of the form $u, u+i, u+2i, \ldots$.

By model  $[a,0,b]$ we mean the situation when the caterpillar spines are given labels $a=a'i$, $b=b'i$,
where $1\le a'<b' <p$. The following lemma plays a key role.

\begin{lemma}\label{lem_a0b} 
We consider the model $[a,0,b]$, $a=a'i$, $b=b'i$, $a'<b'$. Let $\alpha = p-b'-1$, $\beta = b'-a'-1$, $\gamma = a'-1$ such that $\alpha + \beta + \gamma = p-3$. Then:
\begin{enumerate}
    \item  The triple $(\alpha, \beta, \gamma)$  can be realized on the spine component  (see Fig~\ref{a0b}).
\item The triple $(\alpha +1, \beta +1, \gamma +1)$  can be realized on a regular component.

\item If  additionally  $\beta<\gamma<\alpha$,
then on the spine component we can  realize triples $(\alpha -1, \beta +1, \gamma)$ and  $(\alpha -2, \beta +1, \gamma +1)$.

\item Furthermore, if $b'+2a'\ge p+1$ we can realize on the spine component the triple $(\alpha -1, \beta +2, \gamma -1)$.
\end{enumerate}
\end{lemma}

\begin{figure}
\psset{unit=1cm}
\psset{radius=0.2}

\begin{pspicture}(13.5,1.5)

\put(0.5,0.5)
{\begin{pspicture}(0,0)

\put(1,0.3){$\overbrace{\hspace{3cm}}$}
\put(2,0,7){$\gamma ~ \mbox{times}$}
\put(2,0){$\ldots$}

\put(6,0.3){$\overbrace{\hspace{1cm}}$}
\put(6,0.7){$\beta ~ \mbox{times}$}
\put(6.7,0){$\ldots$}
\put(9,0.3){$\overbrace{\hspace{4cm}}$}
\put(10,0.7){$\alpha~ \mbox{times}$}
\put(11.5,0){$\cdots$}

\rput(0,0){\dianode[linestyle=none,fillstyle=solid,fillcolor=yellow]{C1}{\mbox{0}}}
\rput(5,0){\dianode[linestyle=none,fillstyle=solid,fillcolor=yellow]{C4}{\mbox{a}}}
\rput(8,0){\dianode[linestyle=none,fillstyle=solid,fillcolor=yellow]{C7}{\mbox{b}}}

\rput(1,0){\circlenode[linestyle=none,fillstyle=solid,fillcolor=green]{C2}{\mbox{z}}}
\rput(4,0){\circlenode[linestyle=none,fillstyle=solid,fillcolor=green]{C3}{\mbox{z}}}

\rput(6,0){\circlenode[linestyle=none,fillstyle=solid,fillcolor=green]{C5}{\mbox{y}}}
\rput(7,0){\circlenode[linestyle=none,fillstyle=solid,fillcolor=green]{C6}{\mbox{y}}}

\rput(9,0){\circlenode[linestyle=none,fillstyle=solid,fillcolor=green]{C8}{\mbox{x}}}
\rput(13,0){\circlenode[linestyle=none,fillstyle=solid,fillcolor=green]{C9}{\mbox{x}}}

\ncarc{->}{C1}{C2}
\ncarc{->}{C3}{C4}
\ncarc{->}{C4}{C5}
\ncarc{->}{C6}{C7}
\ncarc{->}{C7}{C8}
\ncarc{->}{C9}{C1}

\end{pspicture}}
\end{pspicture}
\caption {One of the possible assignments for the spine component in the model $[a,0,b]$} 
\label{a0b}
\end{figure}

\noindent{\bf Proof.} 
  
{FA:  $x \stackrel{a}{\longrightarrow} y$, $z \stackrel{b}{\longrightarrow} y$, 
$z \stackrel{b-a}{\longrightarrow} x$.}

\textit{1.} One of the possible implementations is shown in Fig.~\ref{a0b}. 
Note that none of the forbidden assignments hold.  Assume first that a vertex $u$ gets label $x$ and the vertex $u+a$ has the label $y$, a contradiction since $\gamma=a'-1$. If a vertex $u$ gets label $z$ and the vertex $u+b$ has the label $y$, a contradiction since $\beta=b'-a'-1<b'$. Finally, assume that a vertex $u$ gets label $z$ and the vertex $u+(b-a)$ has the label $x$,  a contradiction since $\beta=b'-a'-1.$ 

\textit{2.}  Now consider a regular component containing the element $u$. If in the above-described realization we replace $0$ by $u$, $a$ by $u+a$ and $b$ by $u+b$ and assign them labels $x$, $z$ and $y$ respectively, we obtain a realization of the triple $(\alpha +1, \beta +1, \gamma +1)$. 

\textit{3.} {Observe that the assumption $\beta<\gamma<\alpha$ implies the following inequalities:}

$$b'<2a', p+a'>2b', ~~p>a'+b'.$$

We consider the realization shown in Fig.~\ref{a0b}  and the vertex (element) $a+b$. It is easy to see (from the third of the above inequalities) that it lies between $b$ and $0$ 
(according to the orientation generated by $i$). So it has label $x$. Note that changing this label to $y$ does not cause any forbidden assignment. 

In this way, we obtain the realization of the triple $(\alpha -1, \beta +1, \gamma)$. 

Note that from the third of the above inequalities, we can also conclude that the vertex  $2a$ 
lies between $b$ and $0$. So it has the label $x$. This time we change the label to $z$. Again, 
this change does not cause any forbidden assignment. In particular, notice that $2a + (a-b) = a+b$ and this vertex now has the label $y$. Furthermore, $2a\neq a-b$, which vertex 
cannot get label $z$. This follows from the third inequality. 

In this way, we obtain the realization of the triple $(\alpha -2, \beta +1, \gamma +1)$. 

\textit{4}. We return to the situation when the label of vertex $a+b$ is changed from $x$ to $y$.   
It is easy to see that then also vertex $b+2a$ can get the label $y$.
The inequality $b'+2a'\ge p+1$ implies that vertex $b+2a$ lies between $0$ and $a$. \emph{i.e.} it had label $z$ so far. 
In this way, we get the realization of the triple $(\alpha -1, \beta +2, \gamma -1)$. \cbdo

\trou
We continue the proof in this subsection which, as before, will be divided into two parts.

First, let us consider the case when $\alpha + \beta + \gamma = p-3$.

Then this triple can be realized on the spine component by adopting the model $[a,0,b]$ where
$\alpha = p-b'-1$, $\beta = b'-a'-1$, $\gamma = a'-1$ and this follows from the Lemma~\ref{lem_a0b}. \\

Assume now that $\alpha + \beta + \gamma = 2p-3$. In order to realize such a triple, we will need two components.  

Let us further suppose that all three numbers $\alpha, \beta, \gamma$ are odd. We define the numbers $\alpha', \beta', \gamma'$ in the following way:

$$2\alpha'+1=\alpha, ~~2\beta'+1=\beta, ~~2\gamma' + 1=\gamma.$$ 
{Observe that $\alpha'+\beta'+\gamma'=p-3$. }
Thus triple $(\alpha', \beta', \gamma')$ is realized on the spine component in an appropriately chosen model $[a,0,b]$ as above, while the triple $(\alpha'+1, \beta'+1, \gamma'+1)$
is realized in the same model, but on the regular component (cf. Lemma~\ref{lem_a0b} \textit{2.}). Therefore, the triple $(\alpha, \beta, \gamma)$ is realized.

We proceed similarly when only one of the numbers is odd. In particular, if only $\gamma$ is odd, then on the spine component we realize $(\alpha' -1, \beta' +1, \gamma')$ where the numbers are defined by the equalities
 
$$2\alpha'=\alpha, ~~2\beta'+2=\beta, ~~2\gamma' + 1=\gamma,$$
{since $\alpha'+\beta'+\gamma'=p-3$ and $\beta'<\alpha'<\gamma'$}  (cf. Lemma~\ref{lem_a0b} \textit{3.})

Similarly, if only $\alpha$ is odd, then on the spine component we realize $(\alpha' -2, \beta' +1, \gamma' +1)$ where the numbers are defined by the equalities
 
$$2\alpha' -1=\alpha, ~~2\beta'+2=\beta, ~~2\gamma' + 2=\gamma.$$

In the case where only $\beta$ is odd, we must additionally show that the 
inequality required in \textit{4.} holds. Let the numbers $\alpha', \beta', \gamma'$ be defined by the equalities

$$2\alpha'=\alpha, ~~2\beta'+3=\beta, ~~2\gamma'=\gamma.$$

Note that $\beta +\gamma \ge p-2$ because otherwise, the sum $\alpha +\beta +\gamma$ would be $p-3$ and not $2p-3$. Then also $\alpha \le p-1$. So we have $2\beta'+3 +2\gamma' \ge   p-2$. Hence, taking into account the fact that $\beta < \gamma$, we get $p-2<4\gamma'$. Thus $\gamma' \ge \frac{p-1}{4}$  since  $p$ is odd. 
From the fact that $2\alpha' \le p-1$ we get $\alpha' \le \frac{p-1}{2}$.

Therefore, by Lemma~\ref{lem_a0b}, the triple $(\alpha', \beta', \gamma')$ can be realized in the model $[a,0,b]$, $b=b'i$, $a=a'i$ where the relations between $a',b'$ and the numbers $\alpha$ and $\gamma$ are given in \textit{1.}. They imply that $a'\ge \frac{p+3}{4}$ and $b'\ge \frac{p-1}{2}$. 

So we have $b'+2a'\ge \frac{p-1}{2} + 2\frac{p+3}{4} = p+1$. This concludes the proof of this case.

\subsection{Cases when there is no rainbow coloring}

{If the sum $\alpha + \beta + \gamma$ is $p-3$ then such a triple is always realizable because it suffices to use the Lemma~\ref{lem_a0b}. For $\alpha + \beta + \gamma=3p-3$, there is $\alpha=\beta=\gamma=p-1$ and we are done by Case \ref{xeqy}.}

Therefore, we are restricted to the situation $\alpha + \beta + \gamma = 2p-3$.  Let us first consider the case when $\alpha =0$. If $\beta =p-1$ and $\gamma = p-2$, then as we showed above, it is enough to use the model $[a,0,2a]$.
The case $\beta =p-2$ and $\gamma = p-1$ is the subject of the next lemma. 
This fact has been proven in \cite{JamKin}, but for the sake of completeness, we present it here as well.

\subsubsection{Case $(0,p-2,p-1)$ or  $(p-1,p-2,0)$}

\begin{lemma}\label{x=0}
The triple $(0,p-2,p-1)$ is not realizable.
\end{lemma}

\prf Suppose that the triple $(0,p-2,p-1)$ is realizable and let $f$ be the rainbow coloring of the caterpillar. Without loss of generality, let us assume that we use the model $[a,0,b]$

Then the vertices of the caterpillar contain all the elements of the group $A$, and the edges of the caterpillar contain all the elements of the group $A$ except one. Let us denote this missing element by $\zeta$

We will calculate the sum of all the elements of the group in two ways.

Once, as the sum of the elements on the edges (plus one missing element), and the second time as the sum of the elements on the vertices. We have

$$\sum_{v\in X}(a+f(v)) + a + \sum_{v\in Y}(0+f(v)) + b + \sum_{v\in Z}(b+f(v)) + \zeta =$$

$$=\sum_{v\in X}f(v) +a + \sum_{v\in Y}f(v) + b + \sum_{v\in Z}f(v) \pmod p$$

Hence 

$$|X|\cdot a + \sum_{v\in X}f(v) + a + |Y|\cdot 0 +\sum_{v\in Y}f(v) + b + |Z|\cdot b + 
\sum_{v\in Z}f(v) + \zeta =$$

$$= \sum_{v\in X}f(v) +a + \sum_{v\in Y}f(v) + b + \sum_{v\in Z}f(v) \pmod p$$

Since $|X|= 0 \pmod p$ and $|Z|= -1 \pmod p$ we get

$$-b + \zeta=0 ~~\mbox{so}~~\zeta=b.$$
But $b$ appears on the edge of the spine connecting vertices labeled $0$ and $b$, 
a contradiction. \cbdo

\subsubsection{Case $(p-1,0,p-2)$ or  $(p-2,0,p-1)$ where $|Y|=0$} \label{y=0}
We have two cases that are not covered by the considerations so far, $|Y|=0$ and  $|Y|=1$
In this subsection, we deal with the case $|Y|=0$.
By Lemma~\ref{sym} we can assume without loss of generality that $\alpha > \gamma$ that is, we are considering a triple $(p-1,0,p-2)$. Since $|Y|=0<p-1$ by Lemma~\ref{krata} we obtain that if the triple $(\alpha,\beta,\gamma)$ is realizable in the model $[a,0,b]$, then $b\in\langle a\rangle $. But then since no vertex can have label $y$, this means that all regular components are labeled uniformly by Lemma~\ref{xyz}.  No matter how many labels $x$ and $z$ we use on the spine component, there will always be the last $p-1$ of them that cannot be assigned anywhere. Therefore, rainbow coloring does not exist in this case.

\subsubsection{Case $(p-1,1,p-3)$ or  $(p-3,1,p-1)$ where $|Y|=1$}

\begin{figure}
\psset{unit=1cm}
\psset{radius=0.2}

\begin{pspicture}(13.5,3.5)

\put(0.5,2){\begin{pspicture}(0,0)

\put(1.5,0.3){$\overbrace{\hspace{3.5cm}}$}
\put(9.5,0.3){$\overbrace{\hspace{3.5cm}}$}
\put(3,0.7){$\gamma'~ \mbox{times}$}
\put(10,0.7){$\alpha'~ \mbox{times}$}

\put(10.5,0){$\cdots$}
\put(3,0){$\cdots$}

\put(0.5,0.2){$ra$}
\put(12.5,-0.5){$-ra$}

\rput(0,0){\dianode[linestyle=none,fillstyle=solid,fillcolor=yellow]{C1}{\mbox{0}}}
\rput(6.5,0){\dianode[linestyle=none,fillstyle=solid,fillcolor=yellow]{C4}{\mbox{a}}}
\rput(8,0){\dianode[linestyle=none,fillstyle=solid,fillcolor=yellow]{C5}{\mbox{b}}}

\rput(1.5,0){\circlenode[linestyle=none,fillstyle=solid,fillcolor=green]{C2}{\mbox{z}}}
\rput(5,0){\circlenode[linestyle=none,fillstyle=solid,fillcolor=green]{C3}{\mbox{z}}}

\rput(9.5,0){\circlenode[linestyle=none,fillstyle=solid,fillcolor=green]{C6}{\mbox{x}}}

\rput(11.5,0){\circlenode[linestyle=none,fillstyle=solid,fillcolor=green]{C7}{\mbox{x}}}
\rput(13,0){\circlenode[linestyle=none,fillstyle=solid,fillcolor=green]{C8}{\mbox{x}}}

\ncarc{->}{C1}{C2}
\ncarc{->}{C3}{C4}
\ncarc{->}{C4}{C5}
\ncarc{->}{C5}{C6}
\ncarc{->}{C7}{C8}

\end{pspicture}
}
\put(0.5,0.5)
{\begin{pspicture}(0,0)
\put(0,0.3){$\overbrace{\hspace{5cm}}$}
\put(8,0.3){$\overbrace{\hspace{5cm}}$}
\put(3,0.7){$\gamma'+1~ \mbox{times}$}
\put(10,0.7){$\alpha'+1~ \mbox{times}$}

\put(10.5,0){$\cdots$}
\put(3,0){$\cdots$}

\put(0.5,0.1){$ra$}

\put(12.2,-0.5){$u-ra$}

\put(-0.1,-0.5){$u$}
\put(6,-0.5){$u+a$}

\rput(0,0){\circlenode[linestyle=none,fillstyle=solid,fillcolor=green]{C1}{\mbox{z}}}
\rput(6.5,0){\circlenode[linestyle=none,fillstyle=solid,fillcolor=red]{C4}{\mbox{y}}}
\rput(8,0){\circlenode[linestyle=none,fillstyle=solid,fillcolor=green]{C5}{\mbox{x}}}

\rput(1.5,0){\circlenode[linestyle=none,fillstyle=solid,fillcolor=green]{C2}{\mbox{z}}}
\rput(5,0){\circlenode[linestyle=none,fillstyle=solid,fillcolor=green]{C3}{\mbox{z}}}
\rput(9.5,0){\circlenode[linestyle=none,fillstyle=solid,fillcolor=green]{C6}{\mbox{x}}}
\rput(11.5,0){\circlenode[linestyle=none,fillstyle=solid,fillcolor=green]{C7}{\mbox{x}}}
\rput(13,0){\circlenode[linestyle=none,fillstyle=solid,fillcolor=green]{C8}{\mbox{x}}}

\ncarc{->}{C1}{C2}
\ncarc{->}{C3}{C4}
\ncarc{->}{C4}{C5}
\ncarc{->}{C5}{C6}
\ncarc{->}{C7}{C8}

\end{pspicture}}
\end{pspicture}
\caption {The only possible assignments for the spine and regular 
component in the model $[a,0,(r+1)a]$ if $|Y|=1$} 
\label{1y}
\end{figure}

Just consider the case $(p-1,1,p-3)$.  If the triple $(\alpha,\beta,\gamma)$ is realizable in the model $[a,0,b]$, then $b\in\langle a\rangle $ by Lemma~\ref{krata} because $|Y|=1<p-1$. Therefore we can assume that $b=ka$ for some $2\le k \le p-1$. Since only one vertex has the label $y$, it must belong to the regular component. Indeed otherwise, all the regular components would be uniform and our triple would not be possible to implement. 
 Let us assume $r=k-1$ and consider the sequence {$S$} of vertices $S=(0,ra,2ra, \ldots, (p-1)ra)$. { Note that since $p$ is a prime, then $S$ contains all the vertices of the spine component (see Fig~\ref{1y})}. At some position in this sequence there is the element $a$ preceding the element $b$,
(since $b=a+ra$).  The last element of this sequence is the vertex $-ra$.
From the conditions describing forbidden assignments, we know in particular that
$x=ra$ and $z \stackrel{ra}{\longrightarrow} x$ are forbidden. This means that the first element of our sequence after $0$ has the label $z$, and that a vertex with label $z$ cannot precede a vertex with label $x$. 

As a consequence, all vertices {in the sequence $S$} between $0$ and $a$ have the label $z$. 

If any of the vertices lying between $b$ and $-ra$ in $S$ had the label $z$, then the next ones would have to 
have the same label. But the last one is the vertex $-ra$, for which the label $z$ is forbidden. 
Therefore, all vertices between $b$ and $-ra$ (together with him) have the label $x$.

Now consider a regular component containing a unique vertex with a label $y$.
Let $u$ be an arbitrary vertex of this component. We arrange all the vertices in the sequence $S'=(u, u+ra, u+2ra, \ldots,u+(p-1)ra)$. Without loss of generality, we can assume that vertex $u+a$ has label $y$. Therefore the vertex $u$ has the label $z$ (because $x \stackrel{a}{\longrightarrow} y$ is forbidden.
Then similarly, as above, we conclude that all the next vertices in {$S'$}, up to $u+a$ have the label $z$. In turn, the vertices in $S'$ between $u+a$ and the last vertex $u-ra$ have the label $x$. Otherwise, similarly as above, we would conclude that $u-ra$ has label $z$. But then we would have a contradiction because $z \stackrel{b}{\longrightarrow} y$ is forbidden and $u-ra + b =u+a$.

Comparing the numbers of labels $x$ and $z$ on the spine component and the regular component, we see that if on the spine component, we have $\alpha'$ labels $x$ and $\gamma'$ labels $z$ 
($\alpha'+\gamma' = p-3$), then on the regular component there are $\alpha'+1$ and $\gamma'+1$ labels, respectively.

So both $\alpha =2\alpha'+1$ and $\gamma=2\gamma'+1$ are odd numbers, 
so neither of them is equal to either $p-1$ or $p-3$.

\subsection{Cases $p=2$ and $p=3$}

\subsubsection{Case $p=3$}

\begin{theorem}\label{p=3}
Let  $A = {\zet}_3^k$, $k>1$ and let $T$ be a caterpillar with the spine of length two $T\cong C[a_1,a_2, a_3;X,Y,Z]$. 
Then $T$ is rainbow except for the following situations:

1.  $|X|= 0 \pmod 3$ and $|Z| = 2 \pmod 3$ or symmetrically
$|X|= 2 \pmod 3$ and $|Z| = 0 \pmod 3$ 
 
1.  $|Y|=0$ and $|X|= 1 \pmod 3$ and $|Y| = 2 \pmod 3$ 
or symmetrically $|X|= 2 \pmod 3$ and $|Z| = 0 \pmod 3$ 

\end{theorem}

\prf 
Let's first consider the model $[a,0,b]$ where $b\in \langle a\rangle$. Then $b$ must be equal to $-a$. Let's also note that, in this model, all the vertices of the spine component already have labels.

Only regular components are available. Similarly to the case of $p\ge 5$, only triples where 
$\alpha =\gamma$ can be realized on such a component. Then we also have $\alpha=\beta$.

To realize triples where $\alpha \neq\gamma$ we must use the model $[a,0.b]$ 
where $b\notin \langle a\rangle$.

Let's start with remarks about regular components.
We can of course realize uniform sequences on them, i.e. those consisting of only $x$, $y$ or $z$.
Fig.~\ref{43x3} shows that on regular components we can also realize
triples $(3,3,3)$, $(3,6,0)$, $(0,3,6)$ and $(3,0,6)$. 

By Lemma~\ref{sym} we have to consider three cases when the triples have the form $(0,1,2)$, $(2,0,1)$ and $(1,2,0)$.

\noindent{\bf Case $(\alpha,\beta,\gamma)=(0,1,2)$}.
By Lemma~\ref{x=0} this triple is not realizable.

\noindent{\bf Case $(\alpha,\beta,\gamma)=(2,0,1)$}.

If $|Y|=0$ then, this triple is not realizable (see Subsection~\ref{y=0}).
If $|Y|\ge 3$, then the possibility of realization
is shown in Fig.~\ref{2spine3x3}

\noindent{\bf Case $(\alpha,\beta,\gamma)=(1,2,0)$}.
The possibility of realization
is shown in Fig.~\ref{2spine3x3} \cbdo

\begin{figure}
\psset{unit=1cm}
\psset{radius=0.2}

\begin{pspicture}(13.5,6)

\put(0,0.5)
{\begin{pspicture}(0,0)
\psset{unit=0.875cm}

\rput(0,3){\circlenode[linestyle=none,fillstyle=solid,fillcolor=green]{C11}{\mbox{y}}}
\rput(1.5,3){\circlenode[linestyle=none,fillstyle=solid,fillcolor=green]{C12}{\mbox{x}}}
\rput(3,3){\circlenode[linestyle=none,fillstyle=solid,fillcolor=green]{C13}{\mbox{z}}}

\rput(0,1.5){\circlenode[linestyle=none,fillstyle=solid,fillcolor=green]{C21}{\mbox{x}}}
\rput(1.5,1.5){\circlenode[linestyle=none,fillstyle=solid,fillcolor=green]{C22}{\mbox{y}}}
\rput(3,1.5){\circlenode[linestyle=none,fillstyle=solid,fillcolor=green]{C23}{\mbox{z}}}

\rput(0,0){\circlenode[linestyle=none,fillstyle=solid,fillcolor=green]{C31}{\mbox{z}}}
\rput(1.5,0){\circlenode[linestyle=none,fillstyle=solid,fillcolor=green]{C32}{\mbox{y}}}
\rput(3,0){\circlenode[linestyle=none,fillstyle=solid,fillcolor=green]{C33}{\mbox{x}}}

\ncarc{->}{C11}{C12}
\ncarc{->}{C12}{C13}
\ncarc{->}{C13}{C11}

\ncarc{->}{C21}{C22}
\ncarc{->}{C22}{C23}
\ncarc{->}{C23}{C21}

\ncarc{->}{C31}{C32}
\ncarc{->}{C32}{C33}
\ncarc{->}{C33}{C31}

\ncarc{->}{C11}{C21}
\ncarc{->}{C21}{C31}
\ncarc{->}{C31}{C11}

\ncarc{->}{C12}{C22}
\ncarc{->}{C22}{C32}
\ncarc{->}{C32}{C12}

\ncarc{->}{C13}{C23}
\ncarc{->}{C23}{C33}
\ncarc{->}{C33}{C13}

\ncarc{->}{C14}{C24}
\ncarc{->}{C24}{C34}
\ncarc{->}{C34}{C44}
\ncarc{->}{C44}{C54}
\ncarc{->}{C54}{C14}

\end{pspicture}}

\put(3.625,0.5)
{\begin{pspicture}(0,0)
\psset{unit=0.875cm}

\rput(0,3){\circlenode[linestyle=none,fillstyle=solid,fillcolor=green]{C11}{\mbox{x}}}
\rput(1.5,3){\circlenode[linestyle=none,fillstyle=solid,fillcolor=green]{C12}{\mbox{x}}}
\rput(3,3){\circlenode[linestyle=none,fillstyle=solid,fillcolor=green]{C13}{\mbox{x}}}

\rput(0,1.5){\circlenode[linestyle=none,fillstyle=solid,fillcolor=green]{C21}{\mbox{x}}}
\rput(1.5,1.5){\circlenode[linestyle=none,fillstyle=solid,fillcolor=green]{C22}{\mbox{x}}}
\rput(3,1.5){\circlenode[linestyle=none,fillstyle=solid,fillcolor=green]{C23}{\mbox{x}}}

\rput(0,0){\circlenode[linestyle=none,fillstyle=solid,fillcolor=green]{C31}{\mbox{y}}}
\rput(1.5,0){\circlenode[linestyle=none,fillstyle=solid,fillcolor=green]{C32}{\mbox{y}}}
\rput(3,0){\circlenode[linestyle=none,fillstyle=solid,fillcolor=green]{C33}{\mbox{y}}}

\ncarc{->}{C11}{C12}
\ncarc{->}{C12}{C13}
\ncarc{->}{C13}{C11}

\ncarc{->}{C21}{C22}
\ncarc{->}{C22}{C23}
\ncarc{->}{C23}{C21}

\ncarc{->}{C31}{C32}
\ncarc{->}{C32}{C33}
\ncarc{->}{C33}{C31}

\ncarc{->}{C11}{C21}
\ncarc{->}{C21}{C31}
\ncarc{->}{C31}{C11}

\ncarc{->}{C12}{C22}
\ncarc{->}{C22}{C32}
\ncarc{->}{C32}{C12}

\ncarc{->}{C13}{C23}
\ncarc{->}{C23}{C33}
\ncarc{->}{C33}{C13}

\ncarc{->}{C14}{C24}
\ncarc{->}{C24}{C34}
\ncarc{->}{C34}{C44}
\ncarc{->}{C44}{C54}
\ncarc{->}{C54}{C14}

\end{pspicture}}

\put(7.25,0.5)
{\begin{pspicture}(0,0)
\psset{unit=0.875cm}

\rput(0,3){\circlenode[linestyle=none,fillstyle=solid,fillcolor=green]{C11}{\mbox{y}}}
\rput(1.5,3){\circlenode[linestyle=none,fillstyle=solid,fillcolor=green]{C12}{\mbox{y}}}
\rput(3,3){\circlenode[linestyle=none,fillstyle=solid,fillcolor=green]{C13}{\mbox{z}}}

\rput(0,1.5){\circlenode[linestyle=none,fillstyle=solid,fillcolor=green]{C21}{\mbox{y}}}
\rput(1.5,1.5){\circlenode[linestyle=none,fillstyle=solid,fillcolor=green]{C22}{\mbox{y}}}
\rput(3,1.5){\circlenode[linestyle=none,fillstyle=solid,fillcolor=green]{C23}{\mbox{z}}}

\rput(0,0){\circlenode[linestyle=none,fillstyle=solid,fillcolor=green]{C31}{\mbox{y}}}
\rput(1.5,0){\circlenode[linestyle=none,fillstyle=solid,fillcolor=green]{C32}{\mbox{y}}}
\rput(3,0){\circlenode[linestyle=none,fillstyle=solid,fillcolor=green]{C33}{\mbox{z}}}

\ncarc{->}{C11}{C12}
\ncarc{->}{C12}{C13}
\ncarc{->}{C13}{C11}

\ncarc{->}{C21}{C22}
\ncarc{->}{C22}{C23}
\ncarc{->}{C23}{C21}

\ncarc{->}{C31}{C32}
\ncarc{->}{C32}{C33}
\ncarc{->}{C33}{C31}

\ncarc{->}{C11}{C21}
\ncarc{->}{C21}{C31}
\ncarc{->}{C31}{C11}

\ncarc{->}{C12}{C22}
\ncarc{->}{C22}{C32}
\ncarc{->}{C32}{C12}

\ncarc{->}{C13}{C23}
\ncarc{->}{C23}{C33}
\ncarc{->}{C33}{C13}

\ncarc{->}{C14}{C24}
\ncarc{->}{C24}{C34}
\ncarc{->}{C34}{C44}
\ncarc{->}{C44}{C54}
\ncarc{->}{C54}{C14}

\end{pspicture}}

\put(10.875,0.5)
{\begin{pspicture}(0,0)
\psset{unit=0.875cm}

\rput(0,3){\circlenode[linestyle=none,fillstyle=solid,fillcolor=green]{C11}{\mbox{z}}}
\rput(1.5,3){\circlenode[linestyle=none,fillstyle=solid,fillcolor=green]{C12}{\mbox{z}}}
\rput(3,3){\circlenode[linestyle=none,fillstyle=solid,fillcolor=green]{C13}{\mbox{x}}}

\rput(0,1.5){\circlenode[linestyle=none,fillstyle=solid,fillcolor=green]{C21}{\mbox{z}}}
\rput(1.5,1.5){\circlenode[linestyle=none,fillstyle=solid,fillcolor=green]{C22}{\mbox{x}}}
\rput(3,1.5){\circlenode[linestyle=none,fillstyle=solid,fillcolor=green]{C23}{\mbox{z}}}

\rput(0,0){\circlenode[linestyle=none,fillstyle=solid,fillcolor=green]{C31}{\mbox{x}}}
\rput(1.5,0){\circlenode[linestyle=none,fillstyle=solid,fillcolor=green]{C32}{\mbox{z}}}
\rput(3,0){\circlenode[linestyle=none,fillstyle=solid,fillcolor=green]{C33}{\mbox{z}}}

\ncarc{->}{C11}{C12}
\ncarc{->}{C12}{C13}
\ncarc{->}{C13}{C11}

\ncarc{->}{C21}{C22}
\ncarc{->}{C22}{C23}
\ncarc{->}{C23}{C21}

\ncarc{->}{C31}{C32}
\ncarc{->}{C32}{C33}
\ncarc{->}{C33}{C31}

\ncarc{->}{C11}{C21}
\ncarc{->}{C21}{C31}
\ncarc{->}{C31}{C11}

\ncarc{->}{C12}{C22}
\ncarc{->}{C22}{C32}
\ncarc{->}{C32}{C12}

\ncarc{->}{C13}{C23}
\ncarc{->}{C23}{C33}
\ncarc{->}{C33}{C13}

\ncarc{->}{C14}{C24}
\ncarc{->}{C24}{C34}
\ncarc{->}{C34}{C44}
\ncarc{->}{C44}{C54}
\ncarc{->}{C54}{C14}

\end{pspicture}}
\end{pspicture}
\caption {Regular component in the model $[a,0,b]$ where $b \notin \langle a\rangle$} 
\label{43x3}
\end{figure}

\begin{figure}
\psset{unit=1cm}
\psset{radius=0.2}

\begin{pspicture}(13.5,6)

\put(2.5,0.5)
{\begin{pspicture}(0,0)

\rput(0,3){\dianode[linestyle=none,fillstyle=solid,fillcolor=yellow]{C11}{\mbox{0}}}
\rput(1.5,3){\dianode[linestyle=none,fillstyle=solid,fillcolor=yellow]{C12}{\mbox{a}}}
\rput(3,3){\circlenode[linestyle=none,fillstyle=solid,fillcolor=green]{C13}{\mbox{z}}}

\rput(0,1.5){\dianode[linestyle=none,fillstyle=solid,fillcolor=yellow]{C21}{\mbox{b}}}
\rput(1.5,1.5){\circlenode[linestyle=none,fillstyle=solid,fillcolor=green]{C22}{\mbox{y}}}
\rput(3,1.5){\circlenode[linestyle=none,fillstyle=solid,fillcolor=green]{C23}{\mbox{z}}}

\rput(0,0){\circlenode[linestyle=none,fillstyle=solid,fillcolor=green]{C31}{\mbox{z}}}
\rput(1.5,0){\circlenode[linestyle=none,fillstyle=solid,fillcolor=green]{C32}{\mbox{y}}}
\rput(3,0){\circlenode[linestyle=none,fillstyle=solid,fillcolor=green]{C33}{\mbox{x}}}

\ncarc{->}{C11}{C12}
\ncarc{->}{C12}{C13}
\ncarc{->}{C13}{C11}

\ncarc{->}{C21}{C22}
\ncarc{->}{C22}{C23}
\ncarc{->}{C23}{C21}

\ncarc{->}{C31}{C32}
\ncarc{->}{C32}{C33}
\ncarc{->}{C33}{C31}

\ncarc{->}{C11}{C21}
\ncarc{->}{C21}{C31}
\ncarc{->}{C31}{C11}

\ncarc{->}{C12}{C22}
\ncarc{->}{C22}{C32}
\ncarc{->}{C32}{C12}

\ncarc{->}{C13}{C23}
\ncarc{->}{C23}{C33}
\ncarc{->}{C33}{C13}

\ncarc{->}{C14}{C24}
\ncarc{->}{C24}{C34}
\ncarc{->}{C34}{C44}
\ncarc{->}{C44}{C54}
\ncarc{->}{C54}{C14}

\end{pspicture}}

\put(8,0.5)
{\begin{pspicture}(0,0)

\rput(0,3){\dianode[linestyle=none,fillstyle=solid,fillcolor=yellow]{C11}{\mbox{0}}}
\rput(1.5,3){\dianode[linestyle=none,fillstyle=solid,fillcolor=yellow]{C12}{\mbox{a}}}
\rput(3,3){\circlenode[linestyle=none,fillstyle=solid,fillcolor=green]{C13}{\mbox{x}}}

\rput(0,1.5){\dianode[linestyle=none,fillstyle=solid,fillcolor=yellow]{C21}{\mbox{b}}}
\rput(1.5,1.5){\circlenode[linestyle=none,fillstyle=solid,fillcolor=green]{C22}{\mbox{y}}}
\rput(3,1.5){\circlenode[linestyle=none,fillstyle=solid,fillcolor=green]{C23}{\mbox{z}}}

\rput(0,0){\circlenode[linestyle=none,fillstyle=solid,fillcolor=green]{C31}{\mbox{y}}}
\rput(1.5,0){\circlenode[linestyle=none,fillstyle=solid,fillcolor=green]{C32}{\mbox{y}}}
\rput(3,0){\circlenode[linestyle=none,fillstyle=solid,fillcolor=green]{C33}{\mbox{z}}}

\ncarc{->}{C11}{C12}
\ncarc{->}{C12}{C13}
\ncarc{->}{C13}{C11}

\ncarc{->}{C21}{C22}
\ncarc{->}{C22}{C23}
\ncarc{->}{C23}{C21}

\ncarc{->}{C31}{C32}
\ncarc{->}{C32}{C33}
\ncarc{->}{C33}{C31}

\ncarc{->}{C11}{C21}
\ncarc{->}{C21}{C31}
\ncarc{->}{C31}{C11}

\ncarc{->}{C12}{C22}
\ncarc{->}{C22}{C32}
\ncarc{->}{C32}{C12}

\ncarc{->}{C13}{C23}
\ncarc{->}{C23}{C33}
\ncarc{->}{C33}{C13}

\ncarc{->}{C14}{C24}
\ncarc{->}{C24}{C34}
\ncarc{->}{C34}{C44}
\ncarc{->}{C44}{C54}
\ncarc{->}{C54}{C14}

\end{pspicture}}
\end{pspicture}
\caption {Spine component in the model $[a,0,b]$ where $b \notin \langle a \rangle$} 
\label{2spine3x3}
\end{figure}

\subsubsection{Case $p=2$}

\begin{theorem}\label{p=2}
Let  $A = {\zet}_2^k$, $k>1$ and let $T$ be a caterpillar with the spine of length two $T\cong C[a_1,a_2, a_3;X,Y,Z]$. 
Then $T$ is rainbow iff $|X|$ and $|Z|$ are even and $|Y|$ is odd.
\end{theorem}

\prf  Of course, we have to use the model $[a,0,b]$ where $a\neq 0$, $b\neq 0$ and $b\neq a$ 
which in our case means that $b \notin \langle a\rangle$. The only unlabeled vertex in the spine component 
is of the form $a+b$. Since in our case $a+b=a-b=b-a$, this vertex must have the label $y$.

It is easy to check that no triple containing odd numbers can be realized on a regular component. 
In turn, the realizability of the triples $(2,2,0)$, $(2,0, 2)$ and $(0,2,2)$ is shown in
Fig.~\ref{42x2}. This completes the proof. \cbdo

\begin{figure}
\psset{unit=1cm}
\psset{radius=0.2}

\begin{pspicture}(12,2)

\put(0.5,0.5)
{\begin{pspicture}(0,0)
\dotnode(0,0){X1} 
\dotnode(2,0){X2} 
\dotnode(2,2){X3}
\dotnode(0,2){X4}
\ncline{X1}{X2}
\ncline{X2}{X3}
\ncline{X3}{X4}
\ncline{X4}{X1}

\rput(0,0){\dianode[linestyle=none,fillstyle=solid,fillcolor=yellow]{C1}{\mbox{b}}}
\rput(2,0){\circlenode[linestyle=none,fillstyle=solid,fillcolor=green]{C1}{\mbox{y}}}
\rput(0,2){\dianode[linestyle=none,fillstyle=solid,fillcolor=yellow]{C1}{\mbox{0}}}
\rput(2,2){\dianode[linestyle=none,fillstyle=solid,fillcolor=yellow]{C1}{\mbox{a}}}
\end{pspicture}}

\put(4,0.5)
{\begin{pspicture}(0,0)
\dotnode(0,0){X1} 
\dotnode(2,0){X2} 
\dotnode(2,2){X3}
\dotnode(0,2){X4}
\ncline{X1}{X2}
\ncline{X2}{X3}
\ncline{X3}{X4}
\ncline{X4}{X1}

\rput(0,0){\circlenode[linestyle=none,fillstyle=solid,fillcolor=green]{C1}{\mbox{y}}}
\rput(2,0){\circlenode[linestyle=none,fillstyle=solid,fillcolor=green]{C1}{\mbox{y}}}
\rput(0,2){\circlenode[linestyle=none,fillstyle=solid,fillcolor=green]{C1}{\mbox{x}}}
\rput(2,2){\circlenode[linestyle=none,fillstyle=solid,fillcolor=green]{C1}{\mbox{x}}}
\end{pspicture}}

\put(7.5,0.5)
{\begin{pspicture}(0,0)
\dotnode(0,0){X1} 
\dotnode(2,0){X2} 
\dotnode(2,2){X3}
\dotnode(0,2){X4}
\ncline{X1}{X2}
\ncline{X2}{X3}
\ncline{X3}{X4}
\ncline{X4}{X1}

\rput(0,0){\circlenode[linestyle=none,fillstyle=solid,fillcolor=green]{C1}{\mbox{y}}}
\rput(2,0){\circlenode[linestyle=none,fillstyle=solid,fillcolor=green]{C1}{\mbox{z}}}
\rput(0,2){\circlenode[linestyle=none,fillstyle=solid,fillcolor=green]{C1}{\mbox{y}}}
\rput(2,2){\circlenode[linestyle=none,fillstyle=solid,fillcolor=green]{C1}{\mbox{z}}}
\end{pspicture}}

\put(11,0.5)
{\begin{pspicture}(0,0)
\dotnode(0,0){X1} 
\dotnode(2,0){X2} 
\dotnode(2,2){X3}
\dotnode(0,2){X4}
\ncline{X1}{X2} 
\ncline{X2}{X3}
\ncline{X3}{X4}
\ncline{X4}{X1}

\rput(0,0){\circlenode[linestyle=none,fillstyle=solid,fillcolor=green]{C1}{\mbox{z}}}
\rput(2,0){\circlenode[linestyle=none,fillstyle=solid,fillcolor=green]{C1}{\mbox{x}}}
\rput(0,2){\circlenode[linestyle=none,fillstyle=solid,fillcolor=green]{C1}{\mbox{x}}}
\rput(2,2){\circlenode[linestyle=none,fillstyle=solid,fillcolor=green]{C1}{\mbox{z}}}
\end{pspicture}}

\end{pspicture}
\caption {Model $[a,0,b]$ for $p=2$. Spine component and three regular components} \label{42x2}
\end{figure}
\section{Concluding remarks}
We will end this paper with a small remark. Observe that in the proofs of Cases 2.2 and 2.3 of Theorem~\ref{p5}, the crucial point is not that  $p$ is a prime, but rather that $\exp(A)=p\geq 5$ is odd and the group 
$A$ satisfies $|A|\geq 3p$. Therefore, the arguments used in these cases remain valid in a more general setting, even when $\exp(A)\geq 5$ is odd and $|A|\geq 3\exp(A)$. Thus, we can formulate the following proposition:
\begin{proposition} Let $T\cong C(h_1, h_2,h_3)$ be a caterpillar of order $n$ and $A$ be an Abelian group such that $|A|=n$, $\exp(A)=m\geq 5$ is odd and  $|A|\geq 3\exp(A)$. Let $\alpha=h_1 \pmod m$ and  $\beta=h_2 \pmod m$. If   $\beta\geq \alpha$ then  $T$ has an $A$-rainbow labeling.
\end{proposition}

Thus, we finish with the following conjecture.
\begin{conjecture}Let $T\cong C(h_1, h_2,h_3)$ be a caterpillar of order $n$ and $A$ be an Abelian group such that $|A|=n$ and $\exp(A)=m\geq 5$. Let $\alpha=h_1 \pmod m$,  $\beta=h_2 \pmod m$ and $\gamma=h_3\pmod m$. Then $T$ is $A$-rainbow except for the following cases:

\begin{enumerate}
    \item  $\beta=m-2$, and $\alpha=0$, $\gamma=m-1$, or, symmetrically,
$\alpha=m-1$, $\gamma=0$, 

\item $|Y|=0$, and $\alpha=m-1$, $\gamma=m-2$, or, symmetrically,
$\alpha=m-2$, $\gamma=m-1$,

\item $|Y|=1$, and $\alpha=m-1$, $\gamma=m-3$, or, symmetrically,
$\alpha=m-3$, $\gamma=m-1$.\end{enumerate}

\end{conjecture}


\section{Acknowledgements}
This research turned into supported by the AGH University of Krakow under grant no. 16.16.420.054, funded by the Polish Ministry of Science and Higher Education.

\end{document}